\documentclass[12pt]{article}

\usepackage{latexsym,amsmath,amsthm,amssymb,fullpage}

\usepackage{color}

\usepackage{amsthm}

\usepackage{latexsym,amsmath,amsthm,amssymb,fullpage}

\usepackage[utf8]{inputenc}

\usepackage{amssymb}
\usepackage{epstopdf}
\usepackage{amsmath}
\usepackage{amsthm}
\usepackage{amsfonts}
\usepackage{epstopdf}

\newtheorem{theorem}{Theorem}

\newtheorem{corollary}[theorem]{Corollary}

\theoremstyle{definition}

\theoremstyle{remark}

\newcommand\beq{\begin {equation}}
\newcommand\eeq{\end {equation}}
\newcommand\beqs{\begin {equation*}}
\newcommand\eeqs{\end {equation*}}

\newcommand\R{{\mathbb R}}

\newcommand\E{{\mathbb E}}

\title{A dimension-free reverse logarithmic Sobolev inequality for low-complexity functions in Gaussian space}
\author{Ronen Eldan\thanks{Incumbent of the Elaine Blond career development chair. Supported by a European Research Commission Starting Grant and by an Israel Science Foundation grant no. 715/16}\\Weizmann Insitute, Israel \and Michel Ledoux\\University of Toulouse, France
}

\begin {document}

\maketitle

\begin {abstract}
We discuss new proofs, and new forms, of a reverse logarithmic Sobolev inequality,
with respect to the standard Gaussian measure, for low complexity functions, measured in terms of Gaussian-width. In particular, we provide a dimension-free improvement for a related result given in \cite{E18}.
\end {abstract}

\vskip 9mm

\section {A reverse logarithmic Sobolev inequality} \label {sec.1}

The recent work \cite {E18} has put forward a reverse logarithmic Sobolev inequality,
with respect to the standard Gaussian measure, for low complexity functions
measured in terms of Gaussian-width. To briefly recall this inequality, we take again
the notation from \cite {E18}.

Let $\gamma $ denote the standard Gaussian measure on $\R^n$, and let $\nu$
be the probability measure $d\nu = e^f d\gamma$ where $f : \R^n \to \R $ is twice-differentiable.
Let 
$$ 
 \mathrm {D_{KL}} \big ( \nu \, || \,  \gamma \big ) \, = \, \int_{\R^n} f \, d\nu
\qquad \mbox {and} \qquad
\mathcal {I} (\nu) \, = \, \int_{\R^n} |\nabla f |^2 d\nu 
$$
be respectively the Kullback-Leibler divergence (relative entropy) and 
the Fisher information of $\nu$ with respect to $\gamma$, assumed to be finite in the 
following. In particular $\nu$ has a second moment. The standard
logarithmic Sobolev inequality of L.~Gross (cf.~e.g.~\cite {V03,BGL14}) ensures that
\beq \label {eq.logsob}
 \mathrm {D_{KL}} \big ( \nu \, || \, \gamma \big ) \, \leq \, \frac 12 \, \mathcal {I} (\nu) .
\eeq

Let
$$
\mathcal {D} (\nu) \, = \,  \mathbf {GW} (K) 
    \, = \,  \int_{\R^n}  \sup_{t \in K} \, \langle y, t \rangle \, d\gamma (y)
$$
be the Gaussian-width of the set $ K = \{ \nabla f(x) ; x \in \R^n \} $. The quantity 
$\mathcal {D} (\nu)$ is a measure of the complexity of $\nu$ (rather the gradient-complexity of $f$).
It is assumed there that
$ y \mapsto \sup_{t \in K} \, \langle y, t \rangle $ is integrable with respect to $\gamma$.
The following reverse logarithmic Sobolev inequality for measures of low complexity has been established
in \cite [Theorem~4]{E18}. Set $M = M(\nu) :=  - \inf_{x \in \R^n} \Delta f(x) $, assumed to be finite. Then one has
\beq \label {eq.ronen1}
\frac 12 \, \mathcal {I} (\nu)   \ \leq   \, \mathrm {D_{KL}} ( \nu \, || \, \gamma )
+ \frac 12 M_+ +  \mathcal {D}(\nu)^{2/3} \, \mathcal {I} (\nu)^{1/3},
\eeq  
where $M_+ := \max(M,0)$.

Our first theorem gives the following related bound.
\begin {theorem} \label {thm.reverselogsob}
\textsl{ In the preceding notation,}
$$
\frac 12 \, \mathcal {I} (\nu)  \, \leq \,   \mathrm {D_{KL}} \big ( \nu \, || \, \gamma \big)
   + M +  \mathcal {D}(\nu) .
$$
\end {theorem}

As discussed in \cite {E18}, the inequality is sharp on the extremal functions
$f(x) = \langle \alpha , x \rangle$, $x \in \R^n$, $\alpha \in \R^n$, of the
logarithmic Sobolev inequality which have complexity $M = \mathcal {D}(\nu) = 0$. 

To compare between Theorem \ref{thm.reverselogsob} and the bound \eqref{eq.ronen1}, observe that the latter trivially holds true in the case that $\mathcal {I} (\nu) \leq \mathcal {D}(\nu)$, thus we may generally assume that $\mathcal {D}(\nu) \leq \mathcal {D}(\nu)^{2/3} \, \mathcal {I} (\nu)^{1/3}$. Unlike inequality \eqref{eq.ronen1}, the bound of the theorem has the feature that both sides of the inequality are additive with respect to taking products and in this sense it is \emph{dimension-free}. In Section \ref{sec.2} below, we give a slightly different form which improves on equation \eqref{eq.ronen1} and also essentially improves on Theorem \ref{thm.reverselogsob}.

Such a reverse logarithmic Sobolev inequality is of theoretical interest in the study of approximations of partition functions and of low-complexity Gibbs measures on product spaces (cf.~\cite {E18,A18a}). An analogous definition of low-complexity for Boolean functions was considered in \cite{E18}, where it is shown that a low-complexity condition implies that the measure can be decomposed as a mixture of approximate product measures.

In fact, it was very recently shown (\cite{ELS18}) that if a measure $\nu$ satisfies a reverse logarithmic Sobolev inequality, then it is close, in transportation distance, to a mixture of translated Gaussian measures. The combination of such a result with Theorem \ref{thm.reverselogsob} gives a structure theorem for measures of low-complexity, analogous to the one given in \cite{E18}, but for the Gaussian setting. We formulate this as a corollary.

Recall the quadratic Kantorovich metric $\mathrm {W}_2(\nu ,\gamma )$
between $\nu$ and $\gamma$ defined by
$$
\mathrm {W}_2^2(\nu ,\gamma ) \, = \, \inf \int_{\R^n \times \R^n} |x-y|^2 d\pi (x,y)
$$
where the infimum is taken over all couplings $\pi $ 
on $\R^n \times \R^n$ with respective marginals $\nu $ and $\gamma $. A combination of Theorem \ref{thm.reverselogsob} with \cite[Theorem 5]{ELS18} gives,

\begin{corollary}
\textsl{ In the preceding notation, there exists a probability measure $\mu$ such that}
$$
\mathrm {W}_2^2(\nu ,\gamma \star \mu ) \leq 16 n^{1/3} \left (M(\nu) +  \mathcal {D}(\nu)\right )^{2/3}.
$$
\end{corollary}

In particular, the above corollary gives a meaningful result whenever $M(\nu) +  \mathcal {D}(\nu) = o(n)$. It is also conjectured that a dimension-free analogue of \cite[Theorem 5]{ELS18} should hold true, which, combined with our bound would imply the existence of a probability measure $\mu$ such that
$$
\mathrm {W}_2^2(\nu ,\gamma \star \mu ) \leq C \left (M(\nu) +  \mathcal {D}(\nu)\right ),
$$
where $C>0$ is a universal constant.

The proof of \eqref{eq.ronen1} strongly relies on a construction coming from stochastic control theory, 
of an entropy-optimal coupling of the measure $\nu$ to a Brownian motion.
We will come back to it in Section~\ref {sec.2}. In contrast, our proof of Theorem \ref{thm.reverselogsob} follows a simple and direct approach.

\begin {proof}[Proof of Theorem \ref{thm.reverselogsob}]
By integration by parts with respect to the Gaussian measure $\gamma$,
\beqs \begin {split}
\mathcal {I} (\nu)  \, = \, \int_{\R^n} |\nabla f |^2 d\nu 
	& \, = \, \int_{\R^n} |\nabla f |^2 e^f d\gamma \\
	& \, = \, \int_{\R^n} \langle \nabla (e^f) , \nabla f \rangle \, d\gamma \\
	& \, = \, - \int_{\R^n} e^f \, \mathrm {L} f \, d\gamma
	 \, = \, - \int_{\R^n} \mathrm {L} f \, d\nu \\
\end {split} \eeqs
where $\mathrm {L} = \Delta f - \langle x, \nabla f \rangle $ is the Ornstein-Uhlenbeck
operator. Therefore
\beq \label {eq.fisher}
\mathcal {I} (\nu) 
    \, = \,  - \int_{\R^n} \Delta f \, d\nu + \int_{\R^n} \langle x, \nabla f \rangle \, d\nu 
       \, \leq \, M + \int_{\R^n} \langle x, \nabla f \rangle \, d\nu.
\eeq

Let $\pi$ be a coupling on $\R^n \times \R^n$ with respective marginals $\nu $ and $\gamma $. Then,
\beqs \begin {split}
\int_{\R^n} \langle x, \nabla f(x) \rangle \, d\nu(x)
	& \, =\, \int_{\R^n \times \R^n} \langle x, \nabla f(x) \rangle \, d\pi(x,y) \\
	& \, =\, \int_{\R^n \times \R^n} \langle y, \nabla f(x) \rangle \, d\pi (x,y)
		+  \int_{\R^n \times \R^n} \langle x - y, \nabla f(x) \rangle \, d\pi (x,y). \\
\end {split} \eeqs
Now, on the one hand,
$$
\int_{\R^n \times \R^n} \langle y, \nabla f(x) \rangle \, d\pi (x,y)
	 \, \leq \, \int_{\R^n \times \R^n} \sup_{t \in K} \langle y, t \rangle \, d\pi (x,y)
	 \, = \, \int_{\R^n} \sup_{t \in K} \langle y, t \rangle \, d\gamma (y).
$$
On the other hand, by the standard quadratic inequality,
\beqs \begin {split}
\int_{\R^n \times \R^n} \langle x - y, & \nabla f(x)  \rangle \, d\pi (x,y) \\
	& \, \leq \, \frac 12 \int_{\R^n \times \R^n} \big |\nabla f(x)\big|^2 d\pi (x,y)
			+ \frac 12 \int_{\R^n \times \R^n} |x-y|^2 d\pi (x,y) \\
	& \, = \, \frac 12 \int_{\R^n}  \big |\nabla f(x)\big|^2 d\nu (x)
			+ \frac 12 \int_{\R^n \times \R^n} |x-y|^2 d\pi (x,y). \\
\end {split} \eeqs
Taking the infimum over all couplings $\pi $ with respective marginals $\nu $ and $\gamma$,
it holds true that
$$
\int_{\R^n} \langle x, \nabla f \rangle \, d\nu \, \leq \, 	
		\int_{\R^n} \sup_{t \in K} \langle y, t \rangle \, d\gamma (y)
		+ 	\frac 12 \int_{\R^n}  \big |\nabla f(x)\big|^2 d\nu (x)
		+ \frac 12 \, \mathrm{W}_2^2 (\nu, \gamma).
$$ 
Therefore, together with \eqref {eq.fisher},
\beq \label {eq.inter2}
\frac 12 \, \mathcal {I} (\nu)    \, \leq \, M +  \mathcal {D}(\nu)  
		+ \frac 12 \, \mathrm{W}_2^2 (\nu, \gamma).
\eeq

It remains to recall the quadratic transportation cost inequality by M.~Talagrand
(cf.~\cite {V03,BGL14})
\beq \label {eq.talagrand}
\frac 12 \, \mathrm{W}_2^2 (\nu, \gamma) \, \leq \,  \mathrm {D_{KL}} \big ( \nu \, || \, \gamma \big)
\eeq
and the proof is complete.
\end {proof}

\bigskip

Together with the logarithmic Sobolev inequality \eqref {eq.logsob}
$   \mathrm {D_{KL}}  ( \nu || \gamma ) \leq \frac 12 \, \mathcal {I} (\nu) $,
the step \eqref {eq.inter2} of the preceding proof actually also yields
a reverse transportation cost inequality
\beq \label {eq.reversetransport}
 \mathrm {D_{KL}} \big ( \nu \, || \, \gamma \big) \, \leq \, 
		\frac 12 \, \mathrm{W}_2^2 (\nu, \gamma) +  M + \mathcal {D}(\nu) .
\eeq

\bigskip

Theorem~\ref {thm.reverselogsob} may also be deduced from
a classical integrability result for the supremum of a Gaussian
process. Given a set $K \in \R^n$ such that $x \mapsto \sup_{t\in K}\langle x, t \rangle$
is integrable with respect to $\gamma$, setting
$ Z = Z(x) = \sup_{t \in K} \big [ \langle x, t \rangle - \frac 12 |t|^2 \big]$, $ x \in \R^n$,
it holds true that
\beq \label {eq.vitale}
\int_{\R^n} e^{ Z} d\gamma \, \leq \, e^{\int_{\R^n} \sup_{t \in K} \langle x, t \rangle d\gamma}.
\eeq
This inequality was originally put forward in \cite {T82,Vi96} in the context of concentration
properties of suprema of Gaussian processes. Now, the classical entropic inequality (Gibbs variational principle) expresses that
$$
\int_{\R^n} Z \, d\nu 
     \, \leq \,  \mathrm {D_{KL}} \big ( \nu \, || \, \gamma \big) +  \log \int_{\R^n} e^{ Z} d\gamma .
$$
With $K = \{ \nabla f(x) ; x \in \R^n \}$, it therefore follows that
$$
\int_{\R^n}  \Big [ \langle x, \nabla f \rangle - \frac 12 \, |\nabla f |^2 \Big ] d\nu
  \, \leq \, \int_{\R^n} Z \, d\nu 
   \, \leq \, D_{\mathrm {KL}} \big ( \nu \, || \, \gamma \big)
    + \int_{\R^n} \sup_{t \in K} \langle x, t \rangle \, d\gamma.
$$
Again, together with \eqref {eq.fisher}, this yields the conclusion of the theorem.

At the same time, the integrability inequality \eqref {eq.vitale} may be seen as a consequence
of the transportation cost inequality \eqref {eq.talagrand} and the Kantorovich duality.
The argument actually works for any probability $\mu$ on the Borel sets of $\R^n$
satisfying the transportation cost inequality
\beq \label {eq.talagrand2}
 \frac {1}{2C} \, \mathrm{W}_2^2 (\mu', \mu) \, \leq \,   \mathrm {D_{KL}} \big ( \mu' \, || \, \mu \big)
\eeq
for some constant $C >0$ and every $\mu' < \!\! < \mu$ ($C = 1$ for $\mu = \gamma$).

Namely, the Kantorovich duality (cf.~\cite {V03}) expresses that
$$
 \frac 12 \, \mathrm {W}_2^2(\mu' ,\mu ) 
   \, = \, \sup \bigg [ \int_{\R^n} \varphi \, d\mu'  +  \int_{\R^n}  \psi \, d\mu \bigg]
$$
where the supremum runs over the set of measurable functions
$(\varphi, \psi ) \in \mathrm {L}^1(\mu') \times \mathrm {L}^1(\mu)$ satisfying
\beq \label {eq.phipsi}
\varphi (x) + \psi (y) \, \leq \, \frac 12 \, |x-y|^2
\eeq
for $d\mu'$-almost all $x \in \R^n$ and $d\mu$-almost all $y \in \R^n$.
Given then a couple of functions $(\varphi, \psi)$ satisfying \eqref {eq.phipsi},
the choice in \eqref {eq.talagrand2} of $ \frac {d\mu'}{d\mu} = \frac {e^g}{\int_{\R^n} e^g d\mu}$
where
{$g = \frac 1C [\varphi + \int_{\R^n}  \psi \, d\mu ]$} yields that $\log \int_{\R^n} e^g d\mu \leq 0$, that is
$$
\int_{\R^n} e^{\frac 1C \varphi} d\mu \, \leq \, e^{- \frac 1C \int_{\R^n} \psi d\mu}. 
$$

For every $x, y \in \R^n$, and $t \in K$,
$$
\langle x, t \rangle - \frac 12 \, |t|^2 
		\, = \, \langle y, t \rangle  + \langle x-y, t \rangle - \frac 12 \, |t|^2
		\, \leq \, \langle y, t \rangle  + \frac 12 \, |x-y|^2.
$$
Therefore, if
$\varphi (x) = \sup_{t \in K} \big [ \langle x, t \rangle - \frac 12 |t|^2 \big]$, $x \in \R^n$,
then $\psi (y) = - \sup_{t \in K} \langle y, t \rangle$ is a valid candidate for \eqref {eq.phipsi}.
Hence
$$
\int_{\R^n} e^{ \frac 1C Z} d\mu 
\, \leq \, e^{\frac 1C \int_{\R^n} \sup_{t \in K} \langle x, t \rangle d\mu}
$$
which amounts to \eqref {eq.vitale} when $\mu = \gamma$.

\section {Stochastic calculus and the Föllmer process} \label {sec.2}

As mentioned above, the proof of \eqref {eq.ronen1} developed in \cite {E18} uses 
tools from stochastic control theory, and in particular the so-called Föllmer process \cite {F85}
to achieve an entropy-optimal coupling of the measure $\nu$ to a Brownian motion.
This argument has already been proved useful in the study of various
functional inequalities \cite {B02, L13, EL18}.

To summarize a few facts from \cite {L13,E18}, let ${(B_t)}_{t \geq 0}$ be standard Brownian motion
in $\R^n$ (starting from the origin) adapted to a filtration
${(\mathcal {F}_t)}_{t \geq 0}$. Set $v(t,x) = \nabla \log Z(t,x)$, $ t \in [0,1]$, $x \in \R^n$,
where
$$
Z(t,x) \, = \, \E \big ( [e^f](x + B_{1-t} \big).
$$
The Föllmer process ${(X_t)}_{t \in [0,1]}$ solves the stochastic differential equation
$$
X_0 \, = \, 0, \quad dX_t = dB_t + v_t dt
$$
where $v_t = v(t,X_t)$. Amongst its relevant properties, the random variable $X_1$ has
distribution $\nu$, ${(v_t)}_{t \in [0,1]}$ is a martingale, and
$$
\E \bigg ( \int_0^1 |v_t|^2 dt \bigg) \, = \, 2 \, \mathrm {D_{KL}} \big ( \nu \, || \, \gamma \big).
$$

The arguments developed in \cite {E18} thus make use of these properties towards a proof
of the inequality \eqref {eq.ronen1}. Now, actually, a small variation in the same spirit allows for the
following inequality.

\begin {theorem} \label {thm.ronen}
In the notation of Section~\ref {sec.1}, assume that $\nu $ has a finite second moment. Then
$$
   \int_{\R^n}  |x|^2  d\nu -  \int_{\R^n}  |x|^2  d\gamma 
		 \, \leq \, 2 \,  \mathrm {D_{KL}} \big ( \nu \, || \, \gamma \big) 
		 		+  \mathcal {D}(\nu).
$$
\end {theorem}

\begin {proof}
Note first that by integration by parts
$$
    \int_{\R^n}  |x|^2  d\nu - \int_{\R^n}  |x|^2  d\gamma 
    \, = \, \int_{\R^n}  |x|^2  d\nu  - n
   \, = \,  \int_{\R^n} \langle x, \nabla f \rangle \, d\nu
$$
so that the inequality of the theorem amounts to
\beq \label {eq.ronen2}
 \int_{\R^n} \langle x, \nabla f \rangle \, d\nu 
		 \, \leq \,  2 \, \mathrm {D_{KL}} \big ( \nu \, || \, \gamma \big) 
		 		+  \mathcal {D}(\nu).
\eeq

Recall that $K = \{ \nabla f(x) ; x \in \R^n \}$. Arguing as for the proof of Theorem~\ref {thm.reverselogsob}, for
any coupling $\pi$ with respective marginals $\nu$ and $\gamma$,
\beqs \begin {split}
\int_{\R^n} \langle x, \nabla f \rangle \, d\nu 
	& \, = \, \int_{\R^n \times \R^n} \langle x, \nabla f (x) \rangle \, d\pi (x,y) \\
	& \, = \, \int_{\R^n \times \R^n} \langle x - y, \nabla f (x) \rangle \, d\pi (x,y)
				+  \int_{\R^n \times \R^n} \langle y, \nabla f (x) \rangle \, d\pi (x,y) \\
	& \, \leq \, \int_{\R^n \times \R^n} \langle x - y, \nabla f (x) \rangle \, d\pi (x,y)
			+  \int_{\R^n \times \R^n} \sup_{t \in K} \, \langle y, t \rangle \, d\pi (x,y) \\
	& \,  = \, \int_{\R^n \times \R^n} \langle x - y, \nabla f (x) \rangle \, d\pi (x,y)
			+  \mathcal {D}(\nu) .
\end {split} \eeqs
The inequality \eqref {eq.ronen2} would then follow if for some coupling $\pi$,
$$
 \int_{\R^n \times \R^n} \langle x - y, \nabla f (x) \rangle \, d\pi (x,y)
		\, \leq \, 2 \, \mathrm {D_{KL}} \big ( \nu \, || \, \gamma \big).
$$
But the point is that the Föllmer process actually produces an exact coupling for this identity
to hold. Namely, by the definition and properties of ${(X_t)}_{t \in [0,1]}$, $X_1$ has law
$\nu$, $B_1$ has law $\gamma$ and
$$
d \langle X_t - B_t, v_t \rangle \, = \, |v_t|^2 dt  + \langle X_t - B_t, dv_t \rangle.
$$
Since ${(v_t)}_{t\in [0,1]}$ is a martingale, 
$$
\E \big( \langle X_1 - B_1, v_1 \rangle \big) 
		\, = \, \E \bigg ( \int_0^1 |v_t|^2 dt \bigg) 
		\, = \, 2 \, \mathrm {D_{KL}} \big ( \nu \, || \, \gamma \big)
$$
from which the claim follows since $v_1 = \nabla f (X_1)$.
\end {proof}

For the sake of intuition, let us consider an equivalent form of the bound provided by the theorem. Denote,
$$
\mathrm{H}(\nu) := - \int \log \left (\frac{d \nu}{dx} \right ) d \nu,
$$
the differential entropy of $\nu$. It is straightforward to check that the theorem is equivalent to
$$
\mathrm{H}(\nu) - \mathrm{H}(\gamma) \leq \frac{1}{2} \mathcal{D}(\nu).
$$
Note that, in the special case that $\nu$ has the form $\log \frac{d \nu}{d x} = \sup_{t \in K} \big [ \langle x, t \rangle - \frac 12 |t|^2 \big] + \mathrm{const}$, this bound becomes somewhat similar to the bound \eqref{eq.vitale}.

It remains to connect Theorem~\ref {thm.ronen}, or rather inequality \eqref {eq.ronen2}, to Theorem~\ref {thm.reverselogsob}. By the definition of the Fisher information (cf.~\eqref {eq.fisher})
$$
\mathcal {I} (\nu) 
    \, = \,  - \int_{\R^n} \Delta f \, d\nu + \int_{\R^n} \langle x, \nabla f \rangle \, d\nu,
$$
so that \eqref {eq.ronen2} expresses that
\beq \label {eq.reverselogsob2}
\frac 12 \, \mathcal {I} (\nu)  \, \leq \,   \mathrm {D_{KL}} \big ( \nu \, || \, \gamma \big)
   - \frac 12 \int_{\R^n} \Delta f \, d\nu + \frac 12 \,  \mathcal {D}(\nu) .
\eeq
While presented and established with the quantity $M = -\inf_{x \in \R^n} \Delta f(x)$, 
the proof of Theorem~\ref {thm.reverselogsob} shows in the same way that
\beq \label {eq.reverselogsob1}
\frac 12 \, \mathcal {I} (\nu)  \, \leq \,   \mathrm {D_{KL}} \big ( \nu \, || \, \gamma \big)
   -  \int_{\R^n} \Delta f \, d\nu +  \mathcal {D}(\nu) .
\eeq
Hence, if $ \mathcal {D}(\nu) \geq \int_{\R^n} \Delta f \, d\nu$ (which is likely), the inequality
\eqref {eq.reverselogsob2} improves upon \eqref {eq.reverselogsob1}.
On the other hand, it does not seem possible 
to reach \eqref {eq.reverselogsob2} as simply as \eqref {eq.reverselogsob1},
and in any case, the inequality of Theorem~\ref {thm.ronen}, even up to a constant, may not be
deduced from Theorem~\ref {thm.reverselogsob}.

\vskip 8mm

\font\tenrm =cmr10  {\tenrm

\parskip 0mm

\noindent Department of Mathematics

\noindent  Weizmann Institute of Science, 76100 Rehovot, Israel

\noindent roneneldan@gmail.com

\bigskip

\noindent Institut de Math\'ematiques de Toulouse 

\noindent Universit\'e de Toulouse -- Paul-Sabatier, F-31062 Toulouse, France

\noindent \&  Institut Universitaire de France 

\noindent ledoux@math.univ-toulouse.fr

}

\end {document}